\documentclass[11pt]{article}
\usepackage{a4}
\usepackage{amssymb}
\begin{document}

\newcommand{\Section}[1]{\setcounter{equation}{0}\section{#1}}
\renewcommand{\theequation}{\thesection.\arabic{equation}}

\def\be{\begin{equation}}
\def\ee{\end{equation}}
\def\bea{\begin{eqnarray}}
\def\eea{\end{eqnarray}}
\def\A{{\cal A}}
\def\ve{\varepsilon}
\def\ha{\frac{1}{2}}
\def\hR{\hat R}
\def\bx{\bar x}
\def\by{\bar y}
\def\l{\lambda}
\def\hh{{1 \over 2}}
\def\hq{\frac{1}{q}}
\def\nn{\nonumber\\}
\hspace*{\fill}LMU-TPW 98-3\\
\hspace*{\fill}MPI-PhT 98-16\\
\hspace*{\fill}February 1998\\[3ex]
\begin{center}
{\Large \bf \boldmath{$q$}-deformed Hermite Polynomials in \\ 
              \boldmath{$q$}-Quantum Mechanics }\\[3ex]
Ralf Hinterding, Julius Wess\\[3ex]

{\it Sektion Physik der Ludwig-Maximilians-Universit\"at\\                                                                                    
Theresienstr. 37, D-80333 M\"unchen\\ and\\

Max-Planck-Institut f\"ur Physik\\
(Werner-Heisenberg-Institut)\\
F\"ohringer Ring 6, D-80805 M\"unchen\\}

\end{center}

\vskip2cm
%\centerline{\bf Abstract}
%\bigskip
\abstract{ The $q$-special functions appear naturally in $q$-deformed quantum mechanics 
and both sides profit from this fact. Here we study the relation between the $q$-deformed 
harmonic oscillator and the $q$-Hermite polynomials . We discuss: recursion formula, 
generating function, Christoffel-Darboux identity, orthogonality relations and the 
moment functional.}

%\vspace*{\fill}

\newpage

\Section{Introduction}
In the mathematical literature $q$-special functions have been studied very intensively
\cite{Askey}. It turns out that these functions appear in $q$-quantum mechanics when 
we try to diagonalize selfadjoint operators \cite{4names}. Due to the algebraic nature of $q$-quantum 
mechanics many properties of these systems can be derived from the algebraic structure.
This in turn implies special properties of the respective $q$-special functions and would 
be hard to prove differently. A trivial example is the harmonic oscillator ($q=1$) with its 
creation and annihilation operators on the algebraic side and the Hermite polynomials on
the other.

In this short note we generalize this idea to the $q$-deformed harmonic oscillator 
where $q$-deformed Hermite polynomials appear in the eigenfunctions of the 
Hamiltonian \cite{WL2}. A recursion formula for the $q$-Hermite polynomials follows
directly from the construction of the eigenstates.

In chapter 2 we solve this recursion formula explicitly and present a generating 
function of the $q$-Hermite polynomials. We also show that the Christoffel-Darboux
identity follows from the recursion formula. We were, however, not able to prove the
completeness of the polynomials with the help of this identity. From the work on the 
harmonic oscillator \cite{WL2} we actually suspect that the $q$-Hermite polynomials
are not a complete set of functions.

In chapter 3 we give an explicit representation of the eigenstates of the $q$-harmonic 
oscillator in terms of the eigenstates of the coordinates. This then yields orthogonality
relations for the $q$-Hermite polynomials which we derive in chapter 4. It is interesting
that there are two different measures by which the $q$-Hermite polynomials form an 
orthogonal set of functions.  This again indicates the fact that the polynomials are 
not complete.

In chapter 5 we use the matrix elements of powers of the coordinates to define a moment
functional. As expected there are two different measures for this functional but the 
moment functional is independent of the choice of the measure.

Finally we use this moment functional to define an integral and we give its values 
in terms of the $q$-gamma function.

\Section{\boldmath{$q$}-deformed Hermite polynomials}

In the analysis of the $q$-deformed harmonic oscillator, as it was done
in \cite{WL2}, the following recursion formula for the Hermite-polynomials
occurs:
\begin{equation}
 H_{n+1}^{(q)}(\xi) - q^{-\frac{1}{2}}q^{-2n}2\xi H_{n}^{(q)}(\xi) + 2q^{-2}[n]H_{n-1}^{(q)}(\xi) = 0
\end{equation}
The $q$-number $[n]$ is defined as follows:
\be
   [n] \equiv [n]_{q^{-2}} = \frac{1-q^{-2n}}{1-q^{-2}}
\ee
We define the first two polynomials consistent with $H_{-1} = 0$:
\begin{equation}
      H_{0}^{(q)}(\xi)= 1, \qquad  H_{1}^{(q)}(\xi)= 2q^{-\frac{1}{2}}\xi  
\end{equation}
and obtain the next poynomials:
\bea
             H_{2}^{(q)}(\xi)&=& 4q^{-3} \xi^{2} - 2q^{-2} \nonumber     \\
             H_{3}^{(q)}(\xi)&=& 8q^{-\frac{1}{2}}q^{-7}\xi^{3} - 4q^{-\frac{1}{2}}q^{-2}[3]\xi        \\
             H_{4}^{(q)}(\xi)&=& 16q^{-14}\xi^{4} - 8q^{-5}[3](q^{-4} + 1)\xi^{2} + 4q^{-4}[3]    \nonumber  
\eea
A general expression is:
\be
 H_{n}^{(q)}(\xi) = q^{-\frac{n}{2}}\sum_{k=0}^{<\frac{n}{2}]}q^{-k}\frac{q^{-2{n-2k\choose 2}}2^{n-k}(-1)^{k                                       }[n]!}{([2])^{k}[n-2k]![k]_{q^{-4}}!} \xi^{n-2k}  
\ee
The symbol $<\frac{n}{2}]$ means the largest integer smaller or equal $\frac{n}{2}$. 
In the limit $q\to 1$ we obtain from (2.5) the undeformed Hermite polynomials.
A generating function for these $q$-Hermite polynomials has been found in \cite{ralf}:
\begin{equation}
           E_{q^{-2}}(\xi t)e_{q^{-4}}(t^{2}\frac{q}{2(1-q^{2})}) = \sum_{n=0}^{\infty}\frac{q^{\frac{n}{2}}2^{-n}H_{n}^{(q)}(\xi)}{(q^{-2};q^{-2})_{n}}t^{n}
\end{equation}
The exponential functions are defined as follows:
\bea
        E_{q^{-2}}(t) &=& \sum_{n=0}^\infty \frac{q^{-2{n\choose 2}}}{(q^{-2};q^{-2})_n} t^n      \nn
                                \\
        e_{q^{-4}}(t) &=& \sum_{n=0}^\infty \frac{t^n}{(q^{-4};q^{-4})_n }       \nonumber
\eea

In the classical theory of orthogonal polynomials the Christoffel-Darboux identity \cite{che} 
is derived from a recursion relation.

The deduction of this identity for the $q$-deformed Hermite polynomials 
follows exactly the same steps as in the undeformed case. The result is:
\begin{equation}
    \sum_{m=0}^{n} \frac{H_{m}(\xi_\mu^\sigma)H_{m}(\xi_\nu^\tau)}{2^{m}[m]!} = \frac{q^{\frac{1}{2}}q^{2n}}{2^{n+1}[n]!} \frac{H_{n+1}(\xi_\mu^\sigma)H_{n}(\xi_\nu^\tau) - H_{n+1}(\xi_\nu^\tau)H_{n}(\xi_\mu^\sigma)}{(\xi_\mu^\sigma - \xi_\nu^\tau)}
\end{equation}
For the classical polynomials this identity has been used to prove 
completenes of the polynomials \cite{Leb}. In the deformed case the question of 
completeness is still open.

In the mathematical literature the $q$-Hermite II polynomials $\tilde{h}_n$ have been
studied \cite{Askey}. They are related to our polynomials as follows:
\begin{equation}
      H_{n}^{(q)}(\xi) = \frac{q^{-n^{2}}2^{\frac{n}{2}}}{(1-q^{-2})^{\frac{n}{2}}}\tilde{h}_{n}(x';q^{-2})  
\end{equation}
with the rescalation:
\[
       x^\prime = \sqrt{2(q-q^{-1})}\;\xi
\]

\Section{Eigenstates of the \boldmath{$q$}-deformed harmonic oscillator}

Here we are going to exploit the fact that the Hermite polynomials are
part of the Eigenfunctions of the Hamiltonian of the $q$-deformed harmonic
oscillator \cite{WL2}. This oscillator is realized in the Hilbert space of the
$q$-deformed Heisenberg algebra:
\be
       q^\ha XP - q^{-\ha}PX = iU
\ee
The momentum operator has the following eigenvectors and eigenstates:
\bea
            P | l,\sigma \rangle = \sigma q^l | l, \sigma \rangle      \nn
                 l = -\infty \ldots  \infty,  \qquad   \sigma = \pm 1   \\
                \langle l^\prime, \sigma^\prime | l, \sigma \rangle = \delta_{l' l} \delta_{\sigma'\sigma} \nonumber
\eea
For the coordinates we find:
\bea
      X|\nu, \tau \rangle = -\tau \frac{q^{\nu - \ha}}{q-q^{-1}}|\nu, \tau \rangle   \nn
   \nu = -\infty \ldots \infty ,   \qquad \tau = \pm 1    \\
     \langle \nu^\prime, \tau^\prime | \nu , \tau \rangle = 
                \delta_{\nu^\prime \nu} \delta_{\tau^\prime \tau} \nonumber
\eea
The operator $U$ acts on these states as follows:
\bea
   U |l,\sigma \rangle  &=& |l-1, \sigma \rangle   \qquad    \mbox{momentum} \nn
                \\
    U|\nu, \tau \rangle &=& |\nu +1, \tau \rangle   \qquad   \mbox{coordinates}  \nonumber
\eea
These two systems of eigenfunctions are related by the $q$-Fourier 
transformation:
\bea 
  &{}& | 2l, \sigma \rangle = \frac{N_q}{2}\sum_{\nu =-\infty \atop \tau = +,-}^{\infty} q^{\nu +l}\Big\lbrace 
                                              \cos_q2(\nu +l ) - i\sigma\tau\sin_q2(\nu +l )U\Big\rbrace|2\nu,\tau\rangle \nn
                   \\
  &{}& |2l+1, \sigma \rangle = U^{-1}|2l, \sigma \rangle   \nonumber
\eea
The $q$-trigonometric functions are:
\bea 
            \cos_q(2\nu )&\equiv  &\cos(q^{2\nu};q^{-4})     \nn 
                                   \\
             \sin_q(2\nu) &\equiv &\sin(q^{2\nu};q^{-4})   \nonumber
\eea   
and
\be
           N_{q} \equiv \frac{(q^{-2};q^{-4})_{\infty}}{(q^{-4};q^{-4})_{\infty}}
\ee
This is in the notation defined in \cite{Koorn}.The $q$-trigonometric functions 
satisfy the completeness and orthogonality relations:
\bea
          \sum_{n=-\infty}^{+\infty}q^{-2n}\cos(q^{-2(k+n)};q^{-4})\cos(q^{-2(l+n)};q^{-4}) =
                                                        \frac{1}{N_{q}^{2}}q^{2l}\delta_{kl}       \nn
                           \\
          \sum_{n=-\infty}^{+\infty}q^{-2n}\sin(q^{-2(k+n)};q^{-4})\sin(q^{-2(l+n)};q^{-4}) =
                                                        \frac{1}{N_{q}^{2}}q^{2l}\delta_{kl}        \nonumber
\eea
The eigenfunctions of the harmonic oscillator have been defined in \cite{WL2}. 
They are degenerate:
\be
   |n\rangle^r = \frac{1}{\sqrt{2^n[n]!}}H^{(q)}_n(X) |0\rangle^r
\ee
with
\[
              n= 0,1,\ldots \infty,  \qquad   r = 0,1     \nonumber   
\]
The polynomials $H^{(q)}_n(X)$ are functions of the coordinate operator $X$. 
The ground state, however, is easy to define in the momentum representation (3.2):
\be
 |0\rangle ^r = \frac{1}{\sqrt{2}} \sum_{l=-\infty \atop \sigma =+,-}^{\infty} (-1)^l \sigma^{l+r} 
                                       q^{-\ha(l^2 +l)}c_0 |l, \sigma \rangle
\ee
We can Fourier transform these ground states to the $X$ basis using 
(3.5). With the definition:
\[
     c_l  =  q^{-\ha(l^2 +l)}c_0
\] 
the Fourier coefficients are:
\bea
          &{}&\langle 2\nu, \tau|0\rangle^0 = \frac{N_{q}}{\sqrt{2}}\sum_{l=-\infty}^{\infty}
                          q^{\nu+l}\Big(c_{2l}\cos_q2(\nu+l) +i\tau c_{2l+1}\sin_q2(\nu+l)\Big)    \nonumber   \\
                               \nonumber   \\
         &{}&\langle 2\nu +1, \tau |0\rangle ^0 = 0
\eea
and:
\bea
        &{}&\langle 2\nu+1,\tau|0\rangle^1 = -\frac{N_{q}}{\sqrt{2}}\sum_{l=-\infty}^{\infty}
                           q^{\nu+l}\Big(c_{2l+1}q\cos_q2(\nu+l+1) +i\tau c_{2l}\sin_q2(\nu+l)\Big)    \nonumber   \\
                                              \nonumber     \\
        &{}&\langle 2\nu, \tau |0\rangle ^1 = 0  
\eea

\Section{Orthogonality relations}
The eigenstates (3.9) are orthogonal.
\bea
        \delta_{nm}\delta_{rr'} &=& {}^{r'}\!\langle n| m\rangle^r  \nn
                &=& \frac{1}{\sqrt{2^{n+m}[n]! [m]!}}\, {}^{r'}\!\langle 0|H^{(q)}_n(X) H^{(q)}_m(X) |0\rangle^r  \nn
                &=& \frac{1}{\sqrt{2^{n+m}[n]! [m]!}} \sum_{\nu=-\infty \atop \tau=+,-}^\infty 
                                         H^{(q)}_n(\xi_{\nu,\tau})H^{(q)}_m(\xi_{\nu,\tau})\,{}^{r'}\!\langle 0|\nu,\tau\rangle
                                                  \langle \nu,\tau|0\rangle^r       
\eea
with
\[
      \xi_{\nu,\tau} = -\tau\frac{q^{\nu -\ha}}{q-q^{-1}}
\]
We see that depending on $r$ only the even or odd integers $\nu$ contribute to the
sum and we obtain two orthogonality relations:
\bea
     2^n[n]!\delta_{nm} &=& \sum_{\nu=-\infty \atop \tau=+,-}^\infty 
                                                H^{(q)}_n(\xi_{2\nu,\tau})H^{(q)}_m(\xi_{2\nu,\tau})|\langle 2\nu,\tau| 
                                                                   0\rangle^{r=0}|^2      \nn
                                        \\
      2^n[n]!\delta_{nm} &=& \sum_{\nu=-\infty \atop \tau=+,-}^\infty 
                                                H^{(q)}_n(\xi_{2\nu+1,\tau})H^{(q)}_m(\xi_{2\nu+1,\tau})|\langle 2\nu +1, \tau|
                                                                    0\rangle^{r=1}|^2      \nonumber
\eea
These are two orthogonality relations for the $q$-deformed Hermite polynomials
with the two measures:
\bea
     \mu^0(2\nu) = |\langle 2\nu,\tau | 0 \rangle^0|^2    \qquad       \mu^0(2\nu +1) =0   \nn
                                    \\
      \mu^1(2\nu +1) = |\langle 2\nu+1, \tau | 0 \rangle^1|^2     \qquad    \mu^1(2\nu) = 0  \nonumber
\eea
These measures are independant of $\tau$. For the $q$-Hermite II polynomials $\tilde{h}_n$ 
of the mathematical literature (cf. (2.9)) the following 
orthogonality relation is given \cite{Askey}:
\be
      \tilde{N}_q \frac{(q;q)_n}{q^{n^2}} \delta_{nm} =
          \sum_{k=-\infty}^\infty \Big[ \tilde{h}_n(q^k;q)\tilde{h}_m(q^k;q) + 
                                                            \tilde{h}_n(-q^k;q)\tilde{h}_m(-q^k;q)\Big] \omega(q^k)q^k
\ee
here $\tilde{N}_q$ is a normalisation constant independent of $n$ and the summation is 
over all numbers. The measure is  given by:
\be
       \omega(q^k) = \frac{1}{(-q^{2k};q^2)_\infty}
\ee

\Section{The moment functional}

The groundstate expectation value of $\xi^n$ can be computed. We proceed as 
follows: First we expand $\xi^n$ in terms of the $q$-Hermite polynomials:
\be
   \xi^n = \sum_{k=0}^n b_k^{(n)} H_k^{(q)}(\xi)
\ee
The exited states of the harmonic oscillator are given in terms of the 
Hermite polynomials (3.9). They are orthogonal to the groundstate.
We conclude:
\be
       {}^r\!\langle 0| \xi^n | 0 \rangle^r = b^{(n)}_0
\ee
This is independent of $r$. With the help of the generating function of the
$q$-Hermite polynomials (2.6) it is possible to calculate the coefficients
$b^{(n)}_0$  explicitly. Inserting the definition of the  $q$-exponentials
into (2.6) we get:
\bea
  \sum_{k=0}^\infty \frac{q^{-2{k\choose 2}}(\xi t)^k}{(q^{-2};q^{-2})_k}&=&\sum_{n=0}^\infty
                         \frac{q^{\frac{n}{2}}2^{-n}H_n^{(q)}(\xi)t^n}{(q^{-2};q^{-2})_n} \sum_{j=0}^\infty 
                         \frac{q^{-4{j\choose 2}}(-q)^jt^{2j}}{(q^{-4};q^{-4})_j2^j(1-q^2)^j}   \nn
                                                           \nn
         &=& \sum_{n=0}^\infty \sum_{j=0}^n \frac{q^{\frac{n-j}{2}}2^{-(n-j)}H_{n-j}^{(q)}(\xi) 
                            t^{n-j}q^{-4{j\choose 2}}(-1)^jq^jt^{2j}}{(q^{-2};q^{-2})_{n-j} 
                                                    (q^{-4};q^{-4})_j2^j(1-q^2)^j}   \nn
                                                            \nn
         &=& \sum_{k=0}^\infty \sum_{j=0}^{<\frac{k}{2}]}\frac{q^{\frac{k}{2}}2^{-(k-2j)}
                                       H_{k-2j}^{(q)}(\xi)q^{-4{j\choose 2}}(-1)^j}{(q^{-2};q^{-2})_{k-2j}
                                                            (q^{-4};q^{-4})_j2^j(1-q^2)^j} t^k
\eea  
For the last step we changed the summation over $n$ to the summation over
$k$ with $n=k-j$.  On both sides are polynomials in $t$. 
Comparing the coefficients yields:
\bea    
          \frac{q^{-2{k\choose 2}}\xi^k}{(q^{-2};q^{-2})_k} &=& \sum_{j=0}^{<\frac{k}{2}]}
                        \frac{q^{\frac{k}{2}}2^{-k+j}q^{-4{j\choose 2}}(-1)^j}
                               {(q^{-2};q^{-2})_{k-2j}(q^{-4};q^{-4})_j(1-q^2)^j}H_{k-2j}^{(q)}(\xi)    \nn
                                             \nn
              \xi^k &=& \sum_{j=0}^{<\frac{k}{2}]}\frac{q^{\frac{k}{2}}2^{-k+j}q^{-4{j\choose 2}}
                                q^{2{k\choose 2}}(-1)^j(q^{-2};q^{-2})_k}
                                      {(q^{-2};q^{-2})_{k-2j}(q^{-4};q^{-4})_j(1-q^2)^j} H_{k-2j}^{(q)}(\xi)
\eea
This is the linear combination (5.1). Putting the different powers of $q$ together 
we finally get:
\bea
          \xi^k &=& \sum_{j=0}^{<\frac{k}{2}]}\frac{q^{\frac{k}{2}}2^{-k+j}q^{-2j(j-1)}
                            q^{k(k-1)}(-1)^j(q^{-2};q^{-2})_k}{(q^{-2};q^{-2})_{k-2j}
                                (q^{-4};q^{-4})_j(-q^2)^j(1-q^{-2})^j} H_{k-2j}^{(q)}(\xi)      \nn
                                                    \nn
                  &=& \sum_{j=0}^{<\frac{k}{2}]}\frac{q^{k^2-2j^2-\frac{k}{2}}
                                       2^{j-k}(q^{-2};q^{-2})_k}{(q^{-2};q^{-2})_{k-2j}
                                              (q^{-4};q^{-4})_j(1-q^{-2})^j} H_{k-2j}^{(q)}(\xi)
\eea
We see that  for odd powers $\xi^{2n+1}$ the coefficients $b_0^{(2n+1)}$ vanish. 
For even powers we get:
\bea
            b_0^{(2n)} &=& \frac{q^{{2n\choose 2}}(q^{-2};q^{-2})_{2n}}
                                         {2^n(q^{-4};q^{-4})_n(1-q^{-2})^n}              \nn
                                   \nn
                               &=& \frac{q^{{2n\choose 2}}}{2^n}[1][3]\ldots [2n-1]   
\eea
We know how to compute the groundstates in a basis where $\xi$ is 
diagonal, the coefficients for the expansion have been calculated and are
explicitly given by eqns (3.11, 3.12). We obtain:
\be
     b^{(2n)}_0 = \sum_{\nu =-\infty \atop \tau=+,-}^\infty (\xi_{\nu,\tau})^{2n} |\langle \nu,\tau | 0\rangle^r|^2
\ee
More explicitely we find for $r=0$:
\bea
     b^{(2n)}_0 =   N_q^2 \sum_{\nu,j,k =-\infty}^\infty (\xi_{2\nu,\tau})^{2n}
                                 q^{2\nu+j+k}\big(c_{2j}c_{2k}\cos_q2(\nu+j)\cos_q2(\nu+k)     \nn
                                      + c_{2j+1}c_{2k+1}\sin_q2(\nu+j)\sin_q2(\nu+k)\big)
\eea
and for $r=1$:
\bea     
       b^{(2n)}_0 =  N_q^2 \sum_{\nu,j,k =-\infty}^\infty (\xi_{2\nu+1,\tau})^{2n}
                                   q^{2\nu+j+k}\big(c_{2j}c_{2k}\sin_q2(\nu+j)\sin_q2(\nu+k)        \nn
                                        + c_{2j+1}c_{2k+1}q^{2}\cos_q2(\nu+j+1)\cos_q2(\nu+k+1)\big)
\eea
The expansion (5.7) can be interpreted as an integral with the measures:
\be
      d_q\mu^r(\xi) = |\langle \nu,\tau | 0 \rangle^r|^2
\ee
These are two different measures, for $r=0$ the measure is different from 
zero only for even values of $\nu$, for $r=1$ only for odd values of $\nu$.

We have obtained the following moment functionals: 
\bea
             {\cal L}[\xi^{2n}] &=&   \int \xi^{2n} \,d_q\mu^r (\xi) = \frac{q^{2n\choose 2}}{2^n}
                                                                       [1][3]\ldots [2n-1]     \\
                                                                \nn
             {\cal L}[\xi^{2n+1}]&=&   \int \xi^{2n+1} \,d_q\mu^r (\xi) = 0 
\eea
Although we have two different measures the calculation shows 
that the moment functional is independent of the specific measure. 
All that enters into the moment functional is the normalisation of the measure.

That is exactly what is stated by Favards theorem \cite{gas}. It  postulates 
the existence of a unique moment functional for any polynomial sequence 
that is given by a three-term recurrence relation without saying anything 
about the measure, not even about a possible uniqueness.  

All classical orthogonal polynomials are orthogonal with respect 
to a unique measure, but for $q$-polynomials this is not the case 
(e.g. the $q$-Laguerre polynomials) \cite{gas}. It seems that by 
a $q$-quantum mechanical argumentation we have found another example.

\subsection{The \boldmath{$q$}-gamma function}
In this section we want to give the result of the last section - the moment functional - 
in terms of a $q$-deformed gamma function: $\Gamma_q(x)$. This function is 
defined by $(0<q<1)$ \cite{gas}: 
\begin{equation}
            \Gamma_q(x) \equiv \frac{(q;q)_\infty}{(q^x;q)_\infty} (1-q)^{1-x}
\end{equation}
In \cite{gas} also the classical limit $q\to 1$ to the undeformed gamma function 
and some of its properties are given. 

With the help of the identity:
\begin{eqnarray}
            \bigg[\frac{n}{2}\bigg]_{q^{-4}}&=& \frac{1-q^{-2n}}{1-q^{-4}} 
                        = \frac{1-q^{-2n}}{1-q^{-2}}\frac{1-q^{-2}}{1-q^{-4}}      \nonumber   \\
                                    \nonumber    \\
                     &=& \frac{[n]_{q^{-2}}}{[2]_{q^{-2}}} = \frac{[n]}{[2]}
\end{eqnarray}
and the functional equation for the $q$-gamma function we find the result:
\begin{eqnarray}
        \Gamma_{q^{-4}}(\frac{2n+1}{2})&=&\bigg[\frac{2n-1}{2}\bigg]_{q^{-4}}
             \bigg[\frac{2n-3}{2}\bigg]_{q^{-4}}\ldots \bigg[\frac{1}{2}\bigg]_{q^{-4}}
                       \Gamma_{q^{-4}}(\frac{1}{2})    \nonumber  \\
                                        \nonumber    \\
           &=& \frac{[2n-1][2n-3]\ldots  [3][1]\Gamma_{q^{-4}}(\frac{1}{2})}{[2]^n}
\end{eqnarray}
Hence the moment functional for the $q$-Hermite polynomials can be expressed 
by the $q$-gamma function:
\bea
   \int \xi^{2n} \,d_q\mu^r(\xi) &=& \frac{q^{2n\choose 2}}{2^n}
              \frac{[2]^n}{\Gamma_{q^{-4}}(\frac{1}{2})} \Gamma_{q^{-4}}(\frac{2n+1}{2})   \nn
                                                          \\
     \int \xi^{2n+1} \,d_q\mu^r(\xi) &=& 0        \nonumber
\eea

\end{document}